\documentclass[preprint,12pt]{elsarticle}
\usepackage{amsthm,amsmath,amssymb}
\usepackage{graphicx}
\usepackage[colorlinks=true,citecolor=black,linkcolor=black,urlcolor=blue]{hyperref}
\usepackage{mathrsfs}

\theoremstyle{plain}
\newtheorem{theorem}{Theorem}[section]

\theoremstyle{definition}

\theoremstyle{remark}
\newtheorem{remark}[theorem]{Remark}

\allowdisplaybreaks

\begin{document}

\title{Tight Toughness and Isolated Toughness for $\{K_2,C_n\}$-factor critical avoidable graph}
\author{Xiaxia Guan$^{b}$ \footnote{E-mail: guanxiaxia@tyut.edu.cn}, Hongxia Ma $^{b,c}$ \footnote{Corresponding author. E-mail: mjing516@xjnu.edu.cn,} Maoqun Wang$^{d}$ \footnote{E-mail: wangmaoqun@ytu.edu.cn}  \\
\small $~^{a}$ Department of Mathematics, Taiyuan University of Technology, \\ \small Taiyuan 030024, PR China \\ \small $~^{b}$ School of Mathematical Sciences, Xiamen University, \\ \small Xiamen 361000, PR China\\ \small $~^{c}$ School of Mathematical Sciences, Xinjiang Normal University, \\ \small Urumqi 830053, PR China\\ \small $~^{d}$ School of Mathematics and Information Sciences, Yantai University, \\ \small Yantai 264005, PR China\\}
\date{}
\begin{abstract}
A spannning subgraph $F$ of $G$ is a $\{K_2,C_n\}$-factor if each component of $F$ is either $K_{2}$ or $C_{n}$. A graph $G$ is called a $(\{K_2,C_n\},n)$-factor critical avoidable graph if $G-X-e$ has a  $\{K_2,C_n\}$-factor for any $S\subseteq V(G)$ with $|X|=n$ and
$e\in E(G-X)$. In this paper, we first obtain a sufficient condition with regard to isolated toughness of a graph $G$ such that $G$ is $\{K_2,C_{n}\}$-factor critical avoidable. In addition, we give a sufficient condition with regard to tight toughness and isolated toughness of a graph $G$ such that $G$ is $\{K_2,C_{2i+1}|i \geqslant 2\}$-factor critical avoidable respectively.
\end{abstract}

\begin{keyword}
Factor\sep Component factor\sep Toughness\sep Isolated Toughness \sep Factor avoidable graph
\MSC 05C70\sep 05C38\sep 90B10
\end{keyword}
\maketitle

\section{Introduction}
\noindent

Given a graph $G$, let $V(G)$ and $E(G)$ denote the sets of vertices and edges of $G$, respectively. For $S\subseteq V(G)$, we use $G-S$ to
denote the subgraph of $G$ induced by the set $V(G)\setminus S$. For $E'\subseteq E(G)$, we use $G-E'$ to denote the subgraph of $G$ by deleting all edges in $E'$. In particular we write $G-e$ for  $G-E'$ if $E'=\{e\}$.
For $X\subseteq V(G)$, let $c(G-X)$ and $i(G-X)$ denote the numbers of components and isolated vertices in $G-X$, respectively.
All graphs considered in this paper are simple, finite and undirected.

Chv\'{a}tal \cite{Chvatal} introduced the \emph{toughness} of a graph $G$, denoted by $t(G)$, as follows: $t(G)=+\infty$ if $G$ is complete; otherwise
\[t(G)=\min \left \{\frac{|S|}{c(G-S)}\ \bigg|\ S\subseteq V(G), c(G-S)\geq 2 \right\}.\]

Yang, Ma and Liu \cite{Yang} introduced the \emph{isolated toughness} of a graph $G$, denoted by $I(G)$, as follows: $I(G)=+ \infty $ if $G$ is complete; otherwise,
\[I(G)=\min \left \{\frac{|S|}{iso(G-S)}\  \bigg | \ S\subseteq V(G), i(G-S)\geq 2 \right\}.\]

Let $\mathcal{F}$ be a family of graphs. Let $G$ be a graph and let $F$ be a spannning subgraph of $G$. We call that $F$ is an $\mathcal{F}$-\emph{factor} if each component of $F$ is an element of $\mathcal{F}$, we also say that $G$ contains an $\mathcal{F}$-factor.
%In particular, if every member in $\mathcal{F}$ is %a path with at least $k$ vertices, we call $G$ %contains a $P_{\geq k}$-factor; if every member in %$\mathcal{F}$ is either a path with $2$ vertices %or a cycle with at least $2i+1$ vertices, we call %$G$ contains a $\{K_2,C_{2i+1}| i\geq 2\}$-factor;

A graph $G$ is called to be $\mathcal{F}$-\textit{factor avoidable} if $G-e$ has a $\mathcal{F}$-factor for any $e\in E(G)$.
A graph $G$ is called to be $(\mathcal{F},n)$-\emph{factor critical avoidable} if $G-S$ is  $\mathcal{F}$-factor avoidable for any $S\subseteq V(G)$ with $|S|=n$, that is, $G-S-e$ contains an $\mathcal{F}$-factor for any $e\in E(G-S)$.

Let $P_{\geq k}$ denote a family of paths with at least $k$ vertices.
Zhou et al. \cite{ZhouBian} otained some conditions on  binding number of a graph $G$ such that $G$ is $(P_{\geq k},n)$-factor critical avoidable
for $k=2,3$. In \cite{ZhouWu}, Zhou et al. gave a condition on sun toughness of a graph $G$ such that $G$ is $(P_{\geq 3},n)$-factor critical avoidable. Recently, Zhou and Liu \cite{ZhouLiu} gave a condition on toughness and isolated toughness of a graph $G$ such that $G$ is $(P_{\geq k},n)$-factor critical avoidable.

Let $\{K_2,C_{n}\}$ denote a family of paths with $2$ vertices and cycles with at least $n$ vertices.
In this paper, we first obtain a sufficient condition with regard to isolated toughness of a graph $G$ such that $G$ is $\{K_2,C_{n}\}$-factor critical avoidable.
\begin{theorem} \label{cn-thm}
	For two nonnegative integers $n$ and $k$, let $G$ be $(n+k+3)$-connected. If its isolated toughness $I(G)>\frac{n+k+4}{k+3}$, then $G$ is a $\{K_2,C_n\}$-factor critical avoidable graph.
\end{theorem}

Note that if some component of a subgraph $F$ has a cycle of the length $2k$, then there is a subgraph $k$'$K_2$ in this component. Hence, it is enough to consider the cycles whose the length is odd. In addition, we give a sufficient condition with regard to tight toughness and isolated toughness of a graph $G$ such that $G$ is $\{K_2,C_{2i+1}|i \geqslant 2\}$-factor critical avoidable, respectively.
\begin{theorem} \label{tn-thm}
For two nonnegative integers $n$ and $k$, let $G$ be $(n+k+3)$-connected. If its toughness $t(G)>\frac{n+k+4}{k+3}$ or its isolated toughness $I(G)>\frac{3(k+3)+n-1}{k+3}$, then $G$ is a $\{K_2,C_{n}|n\geq 5\}$-factor critical avoidable graph.
\end{theorem}

In Section 2, we shall give the proof of Theorem \ref{cn-thm}. In Section 3, we shall prove Theorem \ref{tn-thm}. Some classes of graphs are constructed to demonstrate the tightness for the bounds our obtain in the following two sections respectively.
\section{ Proof of Theorem \ref{cn-thm}}
\noindent

%Let
%$N_G(S)$ denote the set of all neighbours of the %vertices in $S$.

In this section, we give a proof of Theorem \ref{cn-thm}. We start with a result from \cite{Scheinerman}.
\begin{theorem} \cite{Scheinerman} \label{cn-c}
 A graph $G$ admits a $\{K_2,C_n\}$-factor if and only if $i(G-X)\leq |X|$ for any $X\subset V(G)$.
\end{theorem}

%\begin{theorem} \label{cn-thm}
%For two nonnegative integers $n$ and $k$, let $G$ be $(n+k+3)$-connected. If its isolated toughness $I(G)>\frac{n+k+4}{k+3}$, then $G$ is a $\{K_2,C_n\}$-factor critical avoidable graph.
%\end{theorem}

\emph{Proof of of Theorem \ref{cn-thm}.}
It is easy to see that it is true for a complete graph. We next assume that $G$ is not complete. For any $W\subset V(G)$ with $|W|=n$ and any edge $e\in E(G-W)$, let $H=G-W-e$. We next prove that $H$ admits a $\{K_2,C_n\}$-factor. On the contrary, we assume that $H$ does not admit a $\{K_2,C_n\}$-factor. Then \begin{equation}\label{ig}
i(H-X)\geq |X|+1
\end{equation}
for some subset $X\subset V(H)$ by Theorem \ref{cn-c}.

Note that $G$ is $(n+k+3)$-connected. Then $H$ is $(k+2)$-connected.

\textbf{Claim 1.} $|X|\geq k+3$.

Note that $H$ is  $(k+2)$-connected. We have that $|V(H)|\geq k+2 \geq 2$ and $w(H)=1$. If $X=\emptyset$, then $1=w(H)\geq i(H)\geq 1$, that is, $i(H)=1$ and $V(H)=1$, which contradicts that $|V(H)|\geq 2$.

If $X=\{u\}$, then $|X|=1$ and $i(H-\{u\})\geq 2$  from \eqref {ig}, which again contradicts that $H$ is a $(k+2)$-connected graph.

Assume that $2 \leq|X|\leq k+2$. Then $i(H-X)\geq |X|+1\geq 3$ from \eqref {ig}. So,
\begin{equation*}
i(G-W-X)\geq i(G-W-e-X)-2=i(H-X)-2\geq |X|+1-2\geq 1.
\end{equation*}
This implies that there is an isolated vertex $v\in V(G-W-X)$ in the subgraph $G-W-X$, i.e., $d_{G-W-X}(v)=0$. Thus
\begin{equation*}
d_{G}(v)\leq d_{G-W-X}(v)+|W|+|X|\leq 0+n+k+2\leq n+k+2
\end{equation*}
as $|W|=n$ and $|X|\leq k+2$.
This is a contradiction since $G$ is $(n+k+3)$-connected.

Hence, $|X|\geq k+3$, Claim 1 is verified.

\textbf{Claim 2.} $I(G)\leq\frac{n+k+4}{k+3}$.

It is clear that $i(G-W-e-X) \geq i(G-W-X)\geq i(G-W-e-X)-2$ for any edge $e=uv$. We next distinguish three cases below to prove Claim 2.

\textbf{Case 1.} $i(G-W-X)=i(G-W-X-e)$.

In this case, $d_{G-W-X}(u)>1$ and $d_{G-W-X}(v)>1$, that is, neither $u$ nor $v$ is an isolated vertex in $H-X$. We have that
\begin{eqnarray*}
	I(G) &\leq&\frac{W\cup X}{i(G-W-X)}\\
	&=& \frac{n+|X|}{i(G-W-X-e)}\\
	&=&\frac{n+|X|}{i(H-X)}\\
	&\leq& \frac{n+|X|}{|X|+1}\\
	&=&1+\frac{n-1}{|X|+1}
\end{eqnarray*}
from \eqref {ig}.

If $n\leq 1$, i.e., $n-1\leq 0$, then
\begin{equation*}
I(G)\leq 1<\frac{k+4}{k+3}.
\end{equation*}

If $n\geq 2$, then
\begin{equation*}
I(G)\leq 1+\frac{n-1}{|X|+1}\leq 1+\frac{n-1}{k+3+1}=\frac{n+k+3}{k+4}<\frac{n+k+4}{k+3}
\end{equation*}
from Claim 1.

\textbf{Case 2.} $i(G-W-X)=i(G-W-X-e)-1$.

In this case, there is a vertex $v$ with $d_{G-W-X}(v)=1$. Assume that $u$ is the unique vertex adjacent to $v$ in the subgraph $G-W-X$. Then by \eqref {ig} and Claim 1, we have that
\begin{eqnarray*}
i(G-W-X-u) &=&i(G-W-X)+1\\
&=& i(G-W-X-e)-1+1\\
&=&i(H-X)\\
&\geq& |X|+1\\
&\geq& r+4.
\end{eqnarray*}
Thus,
\begin{eqnarray*}
I(G)& \leq & \frac{W\cup X\cup \{u\}}{i(G-W-X-u)} \\
&\leq& \frac{n+|X|+1}{|X|+1}\\
&=&1+\frac{n}{|X|+1}\\
&\leq& 1+\frac{n}{k+4}\\
&=& \frac{n+k+4}{k+4}\\
&<& \frac{n+k+4}{k+3}.
\end{eqnarray*}
\textbf{Case 3.} $i(G-W-X)=i(G-W-X-e)-2$.

In this case, there is a $K_2$  component in the subgraph $G-W-X$ with $e=E(K_2)$. Then by \eqref {ig} and Claim 1, we have that
\begin{eqnarray*}
i(G-W-X-u) &=&i(G-W-X)+1\\
&=& i(G-W-X-e)-2+1\\
&=&i(H-X)-1\\
&\geq& |X|\\
&\geq& r+3.
\end{eqnarray*}
Thus,
\begin{eqnarray*}
I(G)& \leq & \frac{W\cup X\cup \{u\}}{i(G-W-X-u)} \\
&=& \frac{n+|X|+1}{|X|}\\
&=&1+\frac{n}{|X|}\\
&\leq& 1+\frac{n}{r+3}\\
&=& \frac{n+k+4}{k+3}.
\end{eqnarray*}
Hence, Claim 2 is verified. This is a contradiction since  $I(G)>\frac{n+k+4}{k+3}$. This implies that $G$ is $\{K_2,C_n\}$-factor critical avoidable.

\begin{remark}\label{tough}
If $G$ is $(n+k+3)$-connected, then $I(G)>\frac{n+k+4}{k+3}$ in Theorem \ref{cn-thm} is best possible.

For two nonnegative integers $n$ and $k$ with $n\geq k+1$, let $G=K_{n+k+3}+((k+2)K_1 \cup K_2)$. It is easy to see that $G$ is $(n+k+3)$-connected and $I(G)=\frac{n+k+4}{k+3}$. Assume that $W\subset V(K_{n+k+3})\subset V(G)$ with $|W|=n$ and $e=E(K_2)$. Then $G-W-e=K_{k+3}+((k+4)K_1)$. Let $X=V(K_{k+3})$ in the subgraph $G-W-e$. Then $i(G-W-e-X)=k+4>k+3=|X|$. This implies that  $G$ is not a $\{K_2,C_n\}$-factor critical avoidable graph by Theorem \ref{cn-c}.
\end{remark}

\begin{remark}\label{connected}
The condition on $(n+k+3)$-connected of Theorem \ref{cn-thm} is best possible.

For two nonnegative integers $n$ and $k$ with $n\geq k+1$, let $G=K_{n+k+2}+((k+1)K_1 \cup K_2)$. It is easy to see that $G$ is a $(n+k+2)$-connected graph and $I(G)=\frac{n+k+3}{k+2}>\frac{n+k+4}{k+3}$. Assume that $W\subset V(K_{n+k+2})\subset V(G)$ with $|W|=n$ and $e=E(K_2)$. Then $G-W-e=K_{k+2}+((k+3)K_1)$. Let $X=V(K_{k+2})$ in the subgraph $G-W-e$. Then $i(G-W-e-X)=k+3>k+2=|X|$. This implies that  $G$ is not a $\{K_2,C_n\}$-factor critical avoidable graph by Theorem \ref{cn-c}.
\end{remark}
\section{Proof of Theorem \ref{tn-thm}}
\noindent

In this section, we give a proof of Theorem \ref{tn-thm}.

A connected graph $G$ is called \emph{triangular-cactus} if each block of $G$ is a triangle. It is obvious that $K_3$ is a triangular-cactus. We also say that $K_1$ is a triangular-cactus. For a graph $G$, its the number of triangular-cacti components  is denoted by $c_{tc}(G)$. We start with a result from \cite{Cornuejols}.
\begin{theorem} \cite{Cornuejols} \label{ct-c}
 A graph $G$ admits a $\{K_2,C_{2i+1}| i\geq 2\}$-factor if and only if $c_{tc}(G-X)\leq |X|$ for any $X\subset V(G)$.
\end{theorem}

\emph{Proof of Theorem \ref{tn-thm}.}
It is easy to see that it is true for a complete graph. We next assume that $G$ is not complete. For any $W\subset V(G)$ with $|W|=n$ and any edge $e\in E(G-w)$, let $H=G-W-e$. It is enough to prove that $H$ admits a $\{K_2,C_{2i+1}| i\geq 2\}$-factor. On the contrary, we assume that $H$ does not admit a $\{K_2,C_{2i+1}| i\geq 2\}$-factor. Then
\begin{equation}\label{itc1}
c_{tc}(H-X)\geq |X|+1
\end{equation}
for some subset $X\subset V(H)$ by Theorem \ref{ct-c}.

\textbf{Claim 1.}  $|X|\geq k+3$.

It is obvious that $G-W-e$ is not a complete graph. If $X=\emptyset$, then $c_{tc}(G-W-e)\geq 1$. Note that $1=w(G-W-e)\geq c_{tc}(G-W-e)\geq 1$ since $G-W-e$ is $(k+2)\geq 2$-connected. We have that $w(G-W-e)=1$. It is easy to see that $G-W-e$ is a block, otherwise, which contradicts that $G-W-e$ is $2$-connected. So, $G-W-e$ is a triangle, which contradicts that $G-W-e$ is not a complete graph.

If $X=\{u\}$,  then  $|X|=1$ and $c_{tc}(G-W-e-\{u\})\geq 2$. Note that $w(G-W-e-\{u\})\geq c_{tc}(G-W-e-\{u\})\geq 2$.  This contradicts that $G-W-e$ is a $(k+2)$-connected graph.

Assume that $2 \leq|X|\leq k+2$. Then $w(G-w-e-X)\geq c_{tc}(G-w-e-X)\geq |X|+1\geq 3$. So,
\begin{equation*}
w(G-W-X)\geq w(G-W-e-X)-1\geq 2,
\end{equation*}
that is, $G-W-X$ is not connected.
Note that $|W\cup X|=|W|+|X|\leq n+k+2$. This is a a contradiction since $G$ is $(n+k+3)$-connected.

Hence, $|X|\geq k+3$.

\textbf{Claim 2.} $t(G)\leq\frac{n+k+4}{k+3}$.

By \eqref {itc1} and Claim 1, we have that
\begin{eqnarray*}
w(G-W-X) &\geq&w(G-W-X-e)-1\\
&=& w(H-X)-1\\
&\geq&c_{tc}(H-X)-1\\
&\geq& |X|\\
&\geq& r+3.
\end{eqnarray*}
Thus,
\begin{eqnarray*}
t(G)& \leq & \frac{W\cup X}{w(G-W-X)} \\
&=& \frac{n+|X|}{|X|}\\
&=&1+\frac{n}{|X|}\\
&\leq& 1+\frac{n}{r+3}\\
&=& \frac{n+k+4}{k+3}.
\end{eqnarray*}

\textbf{Claim 3.} $I(G)\leq \frac{3(r+3)+n-1}{r+3}$.

Assume that there are $a$ isolated vertices, $b$ $K_3$'s and $c$ triangular-cactus with at least 5 vertices in the subgraph $H-X$. Let $Y$ be the set of two vertices of every $K_3$.  Denote the triangular-cactus with at least 5 vertices in the subgraph $H-X$ by $H_1, H_2,\ldots, H_c$. Then
\begin{equation*}
c_{tc}(H-X)=a+b+c,
\end{equation*} Let $I_i$ be a maximum independent set of $H_i$ for any $i\in \{1,2,\ldots,c\}$.

\textbf{Claim 3.1.} $|I_i|\geq \frac{|V(H_i)|}{3}$ and $|I_i|\geq 2$.

%We prove it by contradiction. Assume It is obvious that $|I_i|\geq 2$.

%Let $F$ be the triangular-cactus with $|V(F)|=15$ and three vertices in same triangle of degree 4 (see Figure ). Let $\mathcal{F}$ be a family obtained from.

%We will prove that $|I_i|= \frac{|V(H_i)|}{3}$ if $H_i\in \mathcal{F}$, and otherwise, $|I_i|< \frac{|V(H_i)|}{3}$. Assume that there are $f_i$ blocks with $f_i\geq 2$. The conclusion will be proved by induction $f_i$. It is easy to see that it is true for $f_i\leq 4$. Now we assume that $f_i\geq 5$.

%Assume that there are $f_i$ blocks with $f_i\geq 2$. Note that each block is a triangle. We have that $|E(H_i)|=3f_i$ and $|V(H_i)|=2f_i+1$. Assume that there are $g_i$ vertices with the degree 2. Then $$6f_i\geq 2g_i+4(2f_i+1-g_i).$$ So, $g_i\geq f_i+2$. It is obvious that if $v_1$ and $v_2$ are vertices of the degree 2 in two distinct triangles, then  $\{v_1,v_2\}$ is a independent set of $H_i$. Note that there are at most 2 vertices of the degree 2 for a triangle in $H_i$. If $g_i\geq \frac{4f_i+2}{3}$, then  $|I_i|\geq \frac{g_i}{2}\geq \frac{2f_i+1}{3}=\frac{|V(H_i)|}{3}$.

\begin{equation*}
a+b+\left(\sum_{i=1}^c|I_i|\right)\geq a+b+2c\geq a+b+c=c_{tc}(H-X)\geq |X|+1,
\end{equation*}

%\begin{equation*}
%|V(H-X)|\geq a+3b+5c,
%\end{equation*}

%In fact, $c_{tc}(G-W-X-e)-2\leq c_{tc}(G-W-X)\leq c_{tc}(G-W-X-e)+1$.  Then

We next distinguish four cases below to prove the conclusion.

\textbf{Case 1.}  $u,v\in V(H_i)$ for some $i\in \{1,2,\ldots,c\}$ or $u\in V(H_i)$ and $v\in V(H_i)$ for some $i,j\in \{1,2,\ldots,c\}$.

In this case, we have that $c\geq 1$.

\textbf{Sub-Case 1.1.} Assume that $c=1$. Then we claim $i(G-W-X-Y-V(H_1)\setminus I_1)\geq k+4$.  If $|I_1|=2$,  then $|V(H_1)|=5$ and there is an $I_1$ with $u,v\notin I_1$.  Thus,
\begin{eqnarray*}
	&&i(G-W-X-Y-V(H_1)\setminus I_1)\\
	&=&i(G-W-e-X-Y-V(H_1)\setminus I_1)-2\\
	&=& a+b+2\\
	&>& |X|+1\\
	&\geq& k+4.
\end{eqnarray*}

If $|I_i|\geq 3$, then
\begin{eqnarray*}
	&&i(G-W-X-Y-V(H_1)\setminus I_1)\\
	&\geq&i(G-W-e-X-Y-V(H_1)\setminus I_1)-2\\
	&=& i(H-X-Y-V(H_1))-2\\
	&\geq&a+b+|I_1|-2\\
	&\geq& a+b+1\\
	&=& a+b+c\\
	&\geq& |X|+1\\
	&\geq& k+4.
\end{eqnarray*} The conclusion is true.

Thus,
\begin{eqnarray*}
	I(G)& \leq & \frac{W\cup X\cup Y \cup V(H_1)\setminus I_1)}{i(G-W-X-Y-V(H_1)\setminus I_1)} \\
	&=& \frac{n+|X|+2b+(|V(H_1)|-|I_1|)}{a+b+|I_1|}\\
	&\leq& \frac{n+|X|+2a+2b+2|I_1|}{a+b+|I_1|}\\
	&=& 2+\frac{n+|X||}{a+b+|I_1|}\\
	&\leq& 2+\frac{n+|X|}{|X|+1}\\
	&=&3+\frac{n-1}{|X|+1}\\
	&<& 3+\frac{n-1}{k+3}\\
	&=& \frac{3(k+3)+n-1}{k+3}.
\end{eqnarray*}

\textbf{Sub-Case 1.2.} If $c\geq 2$, then $2(c-1)\geq c$.

 Without loss of generality, we may assume that $i=c$. Then
\begin{eqnarray*}
	&&i\left(G-W-X-Y-\sum_{i=1}^{c-1} (V(H_i)\setminus I_i)\right)\\
	&=&i\left(G-W-e-X-Y-\sum_{i=1}^{c-1} (V(H_i)\setminus I_i)\right)\\
	&=& i\left(H-X-Y-\sum_{i=1}^{c-1}(V(H_i)\setminus I_i)\right)\\
	&\geq&a+b+\left(\sum_{i=1}^{c-1}|I_i|\right)\\
	&\geq& a+b+2(c-1)\\
	&\geq& a+b+c\\
	&\geq& |X|+1\\
	&\geq& k+4.
\end{eqnarray*}

Thus,
\begin{eqnarray*}
	I(G)& \leq & \frac{W\cup X\cup Y\cup \sum_{i=1}^{c-1} (V(H_i)\setminus I_i)}{i(G-W-X-Y-\sum_{i=1}^{c-1} (V(H_i)\setminus I_i))} \\
	&=& \frac{n+|X|+2b+\sum_{i=1}^{c-1}(|V(H_i)|-|I_i|)}{a+b+\sum_{i=1}^{c-1}|I_i|}\\
	&=&\frac{n+|X|+2b+\sum_{i=1}^{c-1}|V(H_i)|-\sum_{i=1}^{c-1}|I_i|}{a+b+\sum_{i=1}^{c-1}|I_i|}\\
	&\leq& \frac{n+|X|+2a+2b+3\sum_{i=1}^{c-1}|I_i|-\sum_{i=1}^{c-1}|I_i|}{a+b+\sum_{i=1}^{c-1}|I_i|} \\
	&=& \frac{n+|X|+2a+2b+2\sum_{i=1}^{c-1}|I_i|}{a+b+\sum_{i=1}^{c-1}|I_i|} \\
	&=& 2+\frac{n+|X|}{a+b+\sum_{i=1}^{c-1}|I_i|}\\
	&\leq& 2+\frac{n+|X|}{|X|+1}\\
	&=&3+\frac{n-1}{|X|+1}\\
	&<&3+\frac{n-1}{k+3}\\
	&\leq&\frac{3(k+3)+n-1}{k+3}.
\end{eqnarray*}

\textbf{Case 2.}  $u,v\in(aK_1)$ or $u\in(aK_1)$ and $v\in(H_i)$ for some $i\in \{1,2,\ldots,c\}$. In this case, $a\geq 1$ and
\begin{eqnarray*}
	&&i\left(G-W-X-Y-\{u\}-\sum_{i=1}^c (V(H_i)\setminus I_i)\right)\\
	&=&i\left(G-W-e-X-Y-\{u\}-\sum_{i=1}^c (V(H_i)\setminus I_i)\right)\\
	&=& i\left(H-X-Y-\{u\}-\sum_{i=1}^c (V(H_i)\setminus I_i)\right)\\
	&\geq&a-1+b+\left(\sum_{i=1}^c|I_i|\right).
\end{eqnarray*}

\textbf{Sub-Case 2.1.}
If $a=1$, then  $c\geq 1$. We have that
\begin{eqnarray*}
	a-1+b+\left(\sum_{i=1}^c|I_i|\right) &\geq&a-1+b+2c\\
	&=&a+b+c+c-1\\
	&\geq&a+b+c\\
	&\geq& |X|+1\\
	&\geq& k+4.
\end{eqnarray*}

Thus,
\begin{eqnarray*}
	I(G)& \leq & \frac{W\cup X\cup Y\cup \{u\} \cup\sum_{i=1}^c (V(H_i)\setminus I_i)}{i(G-W-X-Y-\{u\}-\sum_{i=1}^c (V(H_i)\setminus I_i))} \\
	&=& \frac{n+|X|+2b+1+\sum_{i=1}^c(|V(H_i)|-|I_i|)}{a-1+b+\sum_{i=1}^c|I_i|}\\
	&\leq&  \frac{n+|X|+1+2(a-1)+2b+2\sum_{i=1}^c|I_i|}{a-1+b+\sum_{i=1}^c|I_i|} \\
	&=& 2+\frac{n+|X|+1}{a-1+b+\sum_{i=1}^c|I_i|}\\
	&\leq& 2+\frac{n+|X|+1}{|X|+1}\\
	&=&3+\frac{n}{|X|+1}\\
	&<& \frac{3(k+3)+n-1}{k+3}.
\end{eqnarray*}

\textbf{Sub-Case 2.2.}  If $a\geq 2$, then $2(a-1)-2\geq 0$ and
\begin{eqnarray*}
	a-1+b+\left(\sum_{i=1}^c|I_i|\right) &\geq&a-1+b+2c\\
	&\geq&a+b+c-1\\
	&\geq& |X|\\
	&\geq& k+3.
\end{eqnarray*}

Thus, similar to Sub-Case 2.1, we have that
\begin{eqnarray*}
	I(G)& \leq & \frac{W\cup X\cup Y\cup \{u\} \cup\sum_{i=1}^c (V(H_i)\setminus I_i)}{i(G-W-X-Y-\{u\}-\sum_{i=1}^c (V(H_i)\setminus I_i))} \\
	&\leq& \frac{n+|X|+2b+1+2\sum_{i=1}^c|I_i|}{a-1+b+\sum_{i=1}^c|I_i|} \\
	&=& \frac{n+|X|+2(a-1)-2+2b+1+2\sum_{i=1}^c|I_i|}{a-1+b+\sum_{i=1}^c|I_i|} \\
	&=& 2+\frac{n+|X|-1}{a-1+b+\sum_{i=1}^c|I_i|}\\
	&\leq& 2+\frac{n+|X|-1}{|X|}\\
	&=&3+\frac{n-1}{|X|}\\
	&\leq& 3+\frac{n-1}{k+3}\\
	&=& \frac{3(k+3)+n-1}{k+3}.
\end{eqnarray*}

\textbf{Case 3.}
 $u,v\in V(Q)$, or $u\in(aK_3)$ and $v\in(bK_1)$, or $u\in(aK_3)$ and $v\in(H_i)$ for some $i\in \{1,2,\ldots,c\}$,  or $u\in(aK_3)$ and $v\in(Q)$. Let $u\in Y$ if $u\in(aK_3)$ and $u,v\in Y$ if $u,v\in(aK_3)$.

\begin{eqnarray*}
	&&i\left(G-W-X-Y-\sum_{i=1}^c (V(H_i)\setminus I_i)\right)\\
	&=&i\left(G-W-e-X-Y-\sum_{i=1}^c (V(H_i)\setminus I_i)\right)\\
	&=& i\left(H-X-Y-\sum_{i=1}^c (V(H_i)\setminus I_i)\right)\\
	&\geq&a+b+\left(\sum_{i=1}^c|I_i|\right)\\
	&\geq& |X|+1\\
	&\geq& k+4.
\end{eqnarray*}
Thus, similar to Sub-Case 1.2, we have that
\begin{eqnarray*}
	I(G)& \leq & \frac{W\cup X\cup Y\cup \sum_{i=1}^c (V(H_i)\setminus I_i)}{i(G-W-X-Y-\sum_{i=1}^c (V(H_i)\setminus I_i))} \\
	&\leq&  2+\frac{n+|X|}{|X|+1}\\
	&=&3+\frac{n-1}{|X|+1}\\
	&\leq& 3+\frac{n-1}{k+4}\\
	&<&  \frac{3(k+3)+n-1}{k+3}.
\end{eqnarray*}

\textbf{Case 4.}
$u\in(aK_1)$ and $v\in(Q_j)$ for some $j\in \{1,2,\ldots,d\}$. In this case, $a\geq 1$

\begin{eqnarray*}
	&&i\left(G-W-X-Y-\{v\}-\sum_{i=1}^c (V(H_i)\setminus I_i)\right)\\
	&=&i\left(G-W-e-X-Y-\{v\}-\sum_{i=1}^c (V(H_i)\setminus I_i)\right)\\
	&=& i\left(H-X-Y-\{v\}-\sum_{i=1}^c (V(H_i)\setminus I_i)\right)\\
	&\geq&a+b+\left(\sum_{i=1}^c|I_i|\right)\\
	&\geq& |X|+1\\
	&\geq& k+4.
\end{eqnarray*}

Similar to Sub-Case 2.1, we have that
\begin{eqnarray*}
	I(G)& \leq & \frac{W\cup X\cup Y\cup \{v\} \cup\sum_{i=1}^c (V(H_i)\setminus I_i)}{i(G-W-X-Y-\{v\}-\sum_{i=1}^c (V(H_i)\setminus I_i))} \\
	&=& \frac{n+|X|+2b+1+\sum_{i=1}^c(|V(H_i)|-|I_i|)}{a+b+\sum_{i=1}^c|I_i|}\\
	&\leq&  \frac{n+|X|+2a-2+2b+1+2\sum_{i=1}^c|I_i|}{a+b+\sum_{i=1}^c|I_i|} \\
	&=& 2+\frac{n+|X|-1}{a+b+\sum_{i=1}^c|I_i|}\\
	&\leq& 2+\frac{n+|X|-1}{|X|+1}\\
	&=&3+\frac{n-2}{|X|+1}\\
	&<& 3+\frac{n-1}{k+3}\\
	&=& \frac{3(k+3)+n-1}{k+3}.
\end{eqnarray*}
Thus, the conclusion holds.

\begin{remark}\label{tough1}
	If $G$ is $(n+k+3)$-connected, then $t(G)>\frac{n+k+4}{k+3}$ in Theorem \ref{tn-thm} is best possible.
	
	For two nonnegative integers $n$ and $k$ with $n\geq k+1$, let $G=K_{n+k+3}+((k+2)K_1 \cup K_2)$ again. It is easy to see that $G$ is $(n+k+3)$-connected and $t(G)=\frac{n+k+4}{k+3}$. Then  $G$ is not a $\{K_2,C_n\}$-factor critical avoidable graph by Remark \ref{tough}.
\end{remark}

 \begin{remark}\label{tough2}
 	If $G$ is $(n+k+3)$-connected, then $I(G)>\frac{3(k+3)+n-1}{k+3}$ in Theorem \ref{tn-thm} is best possible.
 	
 	For two nonnegative integers $n$ and $k$ with $n\geq k+1$, let $G=K_{n+k+3}+((k+2)K_3 \cup K_2)$. It is easy to see that $G$ is $(n+k+3)$-connected and $I(G)=\frac{3(k+3)+n-1}{k+3}$. Assume that $W\subset V(K_{n+k+3})\subset V(G)$ and $e=E(K_2)$. Then $G-W-e=K_{k+3}+((k+2)K_3 \cup 2K_1)$. Let $X=V(K_{k+3})$ in the subgraph $G-W-e$. Then $c_{tc}(G-W-e-X)=k+4>k+3=|X|$. This implies that  $G$ is not a $\{K_2,C_{2i+1}| i\geq 2\}$-factor critical avoidable graph by Theorem \ref{ct-c}.
 \end{remark}

 \begin{remark}\label{connected1}
 	The condition on $(n+k+3)$-connected of Theorem \ref{tn-thm}  is best possible.
 	
 	For two nonnegative integers $n$ and $k$ with $n\geq k+1$, let $G=K_{n+k+2}+((k+1)K_1 \cup K_2)$.  Then  $G$ is $(n+k+2)$-connected but not a $\{K_2,C_n\}$-factor critical avoidable graph by Remark \ref{connected}. This implies that $G$ is not a $\{K_2,C_{2i+1}| i\geq 2\}$-factor critical avoidable graph, either.
 \end{remark}

\section*{Acknowledgements}
\noindent
This work was supported by the National Natural Science Foundation of China (No. 12301455), the Natural Science Foundation of Shandong Province (No. ZR2023QA080).
\section*{References}
\bibliographystyle{model1b-num-names}
\bibliography{<your-bib-database>}

\end{document}